\input amstex
\documentstyle{amsppt}
\magnification=\magstep1
 \hsize 13cm \vsize 18.35cm \pageno=1
\loadbold \loadmsam
    \loadmsbm
    \UseAMSsymbols
\topmatter
\NoRunningHeads
\title New identities involving $q$-Euler polynomials of higher order
\endtitle
\author
  T. Kim AND Y. H. Kim
\endauthor

\abstract  In this paper, we present new generating functions which
are related to $q$-Euler numbers and polynomials of higher order.
From these generating functions, we give new identities involving
$q$-Euler numbers and polynomials of higher order.
\endabstract
\thanks  2000 AMS Subject Classification: 11B68, 11S80
\newline  The present Research has been conducted by the research
Grant of Kwangwoon University in 2010
\endthanks
\endtopmatter

\document

{\bf\centerline {\S 1. Introduction/ Preliminaries}}

 \vskip 15pt
  Let $\Bbb C$ be the complex number field. We assume that $q\in\Bbb C$ with $|q|<1$ and the $q$-number is defined by $[x]_q=\frac{1-q^x}{1-q}$ in this paper. The $q$-factorial is given by $[n]_q!=[n]_q[n-1]_q\cdots[2]_q[1]_q$ and the $q$-binomial formulae are known that
  $$(x:q)_n=\prod_{i=1}^n(1-xq^{i-1}) =\sum_{i=0}^n {\binom{n}{i}}_q q^{\binom{i}{2}}(-x)^i, \text{ (see [3, 14, 15])},$$
  and
  $$\frac{1}{(x:q)_n}=\prod_{i=1}^n \left(\frac{1}{1-xq^{i-1}}\right)=\sum_{i=0}^{\infty}{\binom{n+i-1}{i}}_q x^i, \text{ (see [3, 5, 14, 15])},$$
where ${\binom{n}{i}}_q=\frac{[n]_q!}{[n-i]_q![i]_q!}=\frac{[n]_q[n-1]_q\cdots[n-i+1]_q}{[i]_q!}. $

The Euler polynomials are defined by $\frac{2}{e^t
+1}e^{xt}=\sum_{n=0}^{\infty}E_n(x)\frac{t^n}{n!}$, for $|t|<\pi$.
In the special case $x=0$, $E_n(=E_n(0))$ are called the $n$-th
Euler numbers. In this paper, we consider the $q$-extensions of
Euler numbers and polynomials of higher order. Barnes' multiple
Bernoulli polynomials are also defined  by
$$\frac{t^r}{\prod_{j=1}^r(e^{a_j t}-1)} e^{xt}=\sum_{n=0}^{\infty}B_n(x,r|a_1, \cdots, a_r)\frac{t^n}{n!},
\text{ where $|t|<\max_{1\leq i \leq r}\frac{2\pi}{|a_i| }$, (see [1, 14]).}\tag1$$ In one of an impressive series of papers (see [1, 6, 14]),
Barnes developed the so-called multiple zeta and multiple gamma
function. Let $a_1, \cdots, a_N$ be positive parameters. Then
Barnes' multiple zeta function is defined  by
$$\zeta_N(s,w| a_1, \cdots, a_N)=\sum_{m_1,\cdots, m_N=0}(w+m_1a_1+\cdots+m_Na_N)^{-s}, \text{ (see [1])},$$
 where $\Re(s)>N$, $ \Re(w)>0$. For $m\in\Bbb Z_+$,  we have
 $$\zeta_N(-m, w|a_1, \cdots, a_N)=\frac{(-1)^mm!}{(N+m)!}B_{N+m}(w,N|a_1,\cdots, a_N).$$

In this paper, we consider Barnes' type multiple $q$-Euler numbers
and polynomials. The purpose of this paper is to present new
generating functions which are related to $q$-Euler numbers and
polynomials of higher order. From the Mellin transformation of these
generating functions, we derive the q-extensions of Barnes' type
multiple zeta functions, which interpolate the $q$-Euler polynomials
of higher order at negative integer. Finally, we give new identities
involving $q$-Euler numbers and polynomials of higher order.

\vskip 10pt

{\bf\centerline {\S 2. $q$-Euler numbers and polynomials of higher order}} \vskip 10pt
 In this section, we assume that $q\in \Bbb C$ with $|q|<1$. Let $x, a_1, \ldots, a_r$ be complex numbers with positive real parts.
   Barnes' type multiple Euler polynomials are defined  by
$$\frac{2^r}{\prod_{j=1}^r( e^{a_jt}+1)}e^{xt}=\sum_{n=0}^{\infty}E_n^{(r)}(x|a_1, \ldots, a_r)\frac{t^n}{n!},
\text{ for $|t|<\max_{1\leq i \leq r}\frac{\pi}{|w_i|},$ (see [6]),}\tag2$$ and
$E_n^{(r)}(a_1, \ldots, a_r)(=E_n^{(r)}(0|a_1, \ldots, a_r))$ are
called the $n$-th Barnes' type multiple Euler numbers. First, we
consider the $q$-extension of Euler polynomials. The $q$-Euler
polynomials are defined by
$$ F_q(t,x)=\sum_{n=0}^{\infty} E_{n,q}(x)\frac{t^n}{n!}=[2]_q\sum_{m=0}^{\infty}(-q)^me^{[m+x]_qt}, \text{ (see [8, 11, 13, 14, 15])}.\tag3$$
  From (3), we have
  $$E_{n,q}(x)=\frac{[2]_q}{(1-q)^n}\sum_{l=0}^n\binom{n}{l}\frac{(-1)^lq^{lx}}{(1+q^{l+1})}.$$
  In the special case $x=0$, $E_{n,q}(=E_{n,q}(0))$ are called the $n$-th $q$-Euler numbers.
From (3), we can easily derive the following relation.
$$E_{0,q}=1, \text{ and } q(qE+1)^n+E_{n,q}=0 \text{ if } n\geq 1, \text{ (see [8, 16, 17])}, $$
 where we use the standard  convention about replacing  $E^k$ by $E_{k,q}.$  It is easy to show that
 $$\lim_{q\rightarrow 1}F_q(t,x)=\frac{2}{e^t +1}e^{xt}=\sum_{n=0}^{\infty}E_n(x)\frac{t^n}{n!}, \text{ (see [2, 3, 19-23])},$$
where $E_n(x)$ are the $n$-th Euler polynomials. For $r\in \Bbb N$,
the Euler polynomials of order $r$ is defined by
  $$\left(\frac{2}{e^t+1}\right)^r e^{xt}=\sum_{n=0}^{\infty}E_n^{(r)}(x)\frac{t^n}{n!}, \text{ for $|t|<\pi$}.\tag4$$
Now we consider the $q$-extension of (4).
$$F_q^{(r)}(t,x)=[2]_q^r\sum_{m_1, \ldots, m_r=0}^{\infty}(-q)^{m_1+\cdots+m_r}e^{[m_1+\cdots+m_r+x]_qt}
=\sum_{n=0}^{\infty}E_{n,q}^{(r)}(x)\frac{t^n}{n!},  \tag5$$ where
$E_{n,q}^{(r)}(x)$ are called the $n$-th $q$-Euler polynomials of
order $r$ (see [10-15]). From (5), we can derive
$$E_{n,q}^{(r)}(x)=\frac{[2]_q^r}{(1-q)^n}\sum_{l=0}^n \binom{n}{l}\frac{(-1)^l q^{lx}}{(1+q^{l+1})^r}.\tag6$$
By (5) and (6), we see that
$$F_q^{(r)}(t,x)=[2]_q^r\sum_{m=0}^{\infty}\binom{m+r-1}{m}(-q)^me^{[m+x]_qt}. \tag7$$
Thus, we note that $\lim_{q\rightarrow
1}F_q^{(r)}(t,x)=\left(\frac{2}{e^t+1}\right)^r
e^{xt}=\sum_{n=0}^{\infty}E_n^{(r)}(x)\frac{t^n}{n!}.$ In the
special case $x=0$, $E_{n,q}^{(r)}(=E_{n,q}^{(r)}(0))$ are called
the $n$-th $q$-Euler numbers of order $r$. By (5), (6) and (7), we
obtain the following proposition.

\proclaim {Proposition 1} For $r\in\Bbb N$, let
$$F_q^{(r)}(t,x)=[2]_q^r\sum_{m_1, \ldots, m_r=0}(-q)^{m_1+\cdots+m_r}e^{[m_1+\cdots+m_r+x]_qt}
=\sum_{n=0}^{\infty}E_{n,q}^{(r)}(x)\frac{t^n}{n!}.$$ Then we have
$$E_{n,q}^{(r)}(x)=\frac{[2]_q^r}{(1-q)^n}\sum_{l=0}^n \binom{n}{l} \frac{(-1)^l q^{lx}}{(1+q^{l+1})^r}
=[2]_q^r\sum_{m=0}^{\infty}\binom{m+r-1}{m}(-q)^m[m+x]_q^n.$$
\endproclaim

From the Mellin transformation of $F_q^{(r)}(t,x)$, we can derive
the following equation.
$$\align
\frac{1}{\Gamma(s)}\int_0^\infty F_q^{(r)}(-t,x)t^{s-1} dt &=[2]_q^r
\sum_{m_1, \ldots, m_r=0}^\infty \frac{(-q)^{m_1 + \cdots +
m_r}}{[m_1 + \cdots + m_r+x]_q^s} \\ &=[2]_q^r \sum_{m=0}^\infty
\binom{m+r-1}{m}(-q)^m \frac{1}{[m+x]_q^s},\tag8
\endalign$$
where $s \in \Bbb C$, $x \neq 0, -1, -2, \ldots$. By (8), we can
define the multiple $q$-zeta function related to $q$-Euler
polynomials.

\proclaim {Definition 2} For $s \in \Bbb C$, $x \in \Bbb R$ with $x
\neq 0, -1, -2, \ldots$, we define the multiple $q$-zeta function
related to $q$-Euler polynomials as $$ \zeta_{q, r} (s, x)=[2]_q^r
\sum_{m_1, \ldots, m_r=0}^\infty \frac{(-q)^{m_1 + \cdots +
m_r}}{[m_1 + \cdots + m_r+x]_q^s}.$$
\endproclaim

Note that $\zeta_{q, r} (s, x)$ is a meromorphic function in whole
complex $s$-plane. From (8), we also note that
$$\zeta_{q, r} (s, x)= [2]_q^r
\sum_{m=0}^\infty \binom{m+r-1}{m}(-q)^m \frac{1}{[m+x]_q^s}.$$ By
Laurent series and the Cauchy residue theorem in (5) and (8),  we
see that
$$\zeta_{q} (-n, x)=E_{n,q}^{(n)}(x), \text{  for } n\in \Bbb Z_+.$$
Therefore, we obtain the following theorem.

\proclaim {Theorem 3} For $r \in \Bbb N, n \in \Bbb Z_+$, and $x \in
\Bbb R$ with $x \neq 0, -1, -2, \ldots$, we have
$$\zeta_q(-n, x)=E_{n,q}^{(r)}(x).$$
\endproclaim

Let $\chi$ be the Dirichlet's character with conductor $f \in \Bbb
N$ with $f \equiv 1 \pmod 2$. Then the generalized $q$-Euler
polynomial attached to $\chi$ are considered by
$$F_{q, \chi} (x)=\sum_{n=0}^\infty E_{n, \chi, q}(x) \frac{t^n}{n!}=[2]_q \sum_{m=0}^\infty (-q)^m \chi(m) e^{[m+x]_q t}.$$
From (3) and (9), we have
$$E_{n, \chi, q}(x) = \frac{[2]_q}{[2]_{q^f}}\sum_{a=0}^{f-1}(-q)^a \chi(a)E_{n, q^f}(\frac{x+a}{f}).$$
In the special case $x=0$, $E_{n, \chi, q}= E_{n, \chi, q}(0)$ are
called the $n$-th generated $q$-Euler number attached to $\chi$.

It is known that the generalized Euler polynomials of order $r$ are
defined by
$$(\frac{2 \sum_{a=0}^{f-1}(-1)^a \chi(a)e^{at}}{e^{ft}+1} )^r e^{xt}=\sum_{n=0}^\infty E_{n,
\chi}^{(r)} (x) \frac{t^n}{n!}, \tag10$$ for $|t|< \frac{\pi}{f}$.

We consider the $q$-extension of (10). The generalized $q$-Euler
polynomials of order $r$ attached to $\chi$ are defined by
$$\align
F_{q, \chi}^{(r)}(t,x)&=[2]_q^r \sum_{m_1, \ldots, m_r=0}^\infty
(-q)^{m_1 + \cdots + m_r}( \prod_{i=1}^r \chi(m_i )) e^{[m_1 +
\cdots + m_r+x]_q t} \\ &=\sum_{n=0}^\infty E_{n, \chi, q}^{(r)} (x)
\frac{t^n}{n!}, \text{ (see [14, 15])}.\tag11
\endalign$$
Note that $$\lim_{q \rightarrow 1}F_{q, \chi}^{(r)}(t,x)=(\frac{2
\sum_{a=0}^{f-1}(-1)^a \chi(a)e^{at}}{e^{ft}+1} )^r. $$ By (11), we
easily see that
$$\align
& E_{n, \chi, q}^{(r)} (x) = \frac{[2]_q^r}{(1-q)^n} \sum_{l=0}^n
\binom{n}{l}(-q^x)^l \sum_{a_1, \ldots, a_r=0}^{f-1} (\prod_{j=1}^r
\chi(a_j)) \frac{(-q^{l+1})^{\sum_{i=1}^r a_i}}{(1+q^{(l+1)f})^r} \\
& \quad =[2]_q^r \sum_{m_1, \ldots, m_r=0}^\infty (-q)^{m_1 + \cdots
+ m_r} (\prod_{i=1}^r \chi(m_i))[m_1 +\cdots +m_r +x]_q^n.
\endalign$$
For $s \in \Bbb C$, $x \in \Bbb R$ with $x \neq 0, -1, -2, \ldots$,
we have
$$\align
&\frac{1}{\Gamma(s)}\int_0^\infty F_{q, \chi}^{(r)}(-t,x)t^{s-1} dt\\
& \qquad=[2]_q^r \sum_{m_1, \ldots, m_r=0}^\infty \frac{(-q)^{m_1 +
\cdots + m_r}(\prod_{i=1}^r \chi(m_i))}{[m_1 + \cdots + m_r+x]_q^s}, \text{ (see [15])}.
\tag12
\endalign$$
From (12), we can consider the Dirichlet's type multiple
$q$-$l$-function as follows :

\proclaim {Definition 4} For $s \in \Bbb C$, $x \in \Bbb R$ with $x
\neq 0, -1, -2, \ldots$, we define the Dirichlet's type multiple
$q$-$l$-function as
$$l_q (s, x |\chi) =[2]_q^r \sum_{m_1, \ldots, m_r=0}^\infty \frac{(-q)^{m_1 +
\cdots + m_r}(\prod_{i=1}^r \chi(m_i))}{[m_1 + \cdots +
m_r+x]_q^s}, \text{ (see [15])}.$$
\endproclaim
By Laurent series and the Cauchy residue theorem in (11) and (12),
we obtain the following theorem.

\proclaim {Theorem 5} For $n \in \Bbb Z_+$, we have
$$l_q (-n, x |\chi)=E_{n, \chi, q}^{(r)}(x).$$
\endproclaim

For $h \in \Bbb Z$ and $r \in \Bbb N$, we consider the extended
$r$-ple $q$-Euler polynomials.
$$\align
F_q^{(h,r)}(t,x)&=[2]_q^r \sum_{m_1, \ldots, m_r=0}^\infty
q^{\sum_{j=1}^{r} (h-j+1)m_j}(-1)^{\sum_{j=1}^{r} m_j} e^{[m_1 +
\cdots + m_r+x]_q t} \\ &=\sum_{n=0}^\infty E_{n, q}^{(h, r)} (x)
\frac{t^n}{n!}.\tag13
\endalign$$
Note that $$\lim_{q \rightarrow 1}F_q^{(h, r)}(t,x)=(\frac{2}{e^t
+1})^r e^{xt} = \sum_{n=0}^{\infty} E_n^{(r)}(x) \frac{t^n}{n!}.$$
From (13), we note that
$$\align
E_{n, q}^{(h, r)} (x) & = \frac{[2]_q^r}{(1-q)^n} \sum_{l=0}^n
\binom{n}{l}\frac{ (-q^x)^l}{(-q^{h-r+l+1} :q)_r} \\ & =[2]_q^r
\sum_{m=0}^\infty \binom{m+r-1}{m}_q (-q^{h-r+1})^m [m+x]_q^n.
\tag14
\endalign$$
By (14), we easily see that
$$F_q^{(h,r)}(t,x)=[2]_q^r
\sum_{m=0}^\infty \binom{m+r-1}{m}_q (-q^{h-r+1})^m e^{[m+x]_q t}, \text{ (see [11, 13, 14])}.
\tag15
$$
Using the Mellin transform for $F_q^{(h,r)}(t,x)$, we have
$$\align
&\frac{1}{\Gamma(s)}\int_0^\infty F_{q}^{(r)}(-t,x)t^{s-1} dt\\
& \qquad=[2]_q^r \sum_{m_1, \ldots, m_r=0}^\infty \frac{(-1)^{m_1 +
\cdots + m_r}q^{\sum_{j=1}^{r} (h-j+1)m_j}}{[m_1 + \cdots +
m_r+x]_q^s}, \text{ (see [13, 14, 15])},\tag16
\endalign$$
for $s \in \Bbb C$, $x \in \Bbb R$ with $x \neq 0, -1, -2, \ldots$.
Now we can define the extended $q$-zeta function associated with
$E_{n, q}^{(h, r)} (x)$.

\proclaim {Definition 6} For $s \in \Bbb C$, $x \in \Bbb R$ with $x
\neq 0, -1, -2, \ldots$, we define the (h, q)-zeta function as
$$\zeta_{q, r}^{(h)} (s, x) =[2]_q^r \sum_{m_1 , \ldots, m_r=0}^\infty
\frac{(-1)^{m_1 +\cdots + m_r}q^{\sum_{j=1}^{r} (h-j+1)m_j}}{[m_1 +
\cdots + m_r+x]_q^s}.$$
\endproclaim

Note that $\zeta_{q, r}^{(h)} (s, x)$ is also a meromorphic function
in whole complex $s$-plane. From (16) and (15), we note that
$$\zeta_{q, r}^{(h)} (s, x)= [2]_q^r
\sum_{m=0}^\infty \binom{m+r-1}{m}_q (-q^{h-j+1})^m
\frac{1}{[m+x]_q^s}. \tag17$$ Using the Cauchy residue theorem and
Laurent series in (16),  we obtain the following theorem.

\proclaim {Theorem 7} For $n \in \Bbb Z_+$, we have
$$\zeta_{q, r}^{(h)}(-n, x)=E_{n,q}^{(h, r)}(x).$$
\endproclaim

We consider the extended $r$-ple generalized $q$-Euler polynomials
as follows :
$$\align
& F_{q, \chi}^{(h,r)}(t,x)\\&=[2]_q^r \sum_{m_1, \ldots,
m_r=0}^\infty
q^{\sum_{j=1}^{r} (h-j+1)m_j}(-1)^{\sum_{j=1}^{r} m_j} (\prod_{j=1}^r \chi(m_j))e^{[m_1 + \cdots + m_r +x]_q t} \tag18\\
&=\sum_{n=0}^\infty E_{n, \chi, q}^{(h, r)} (x)\frac{t^n}{n!}.
\endalign$$
By (18), we see that
$$\align
& E_{n, \chi, q}^{(h, r)} (x) = \frac{[2]_q^r}{(1-q)^n}\sum_{a_1,
\ldots, a_r =0}^{f-1} (-1)^{\sum_{j=1}^r a_j}(\prod_{j=1}^r
\chi(a_j))\sum_{l=0}^n \binom{n}{l}\frac{ (-1)^l q^{lx}
q^{(h-j+l+1)a_j}}{(-q^{(h-r+l+1)f} :q^f)_r} \\ & \,
=\frac{[2]_q^r}{[2]_{q^f}^r}[f]_q^n \sum_{a_1, \ldots, a_r =0}^{f-1}
(-1)^{\sum_{j=1}^r a_j}(\prod_{j=1}^r \chi(a_j)) q^{\sum_{j=1}^r
(h-j+1)a_j} \zeta_{q^f, r}^{(h)}(-n, \frac{x+\sum_{j=1}^r a_j}{f}).
\tag19
\endalign$$
Therefore, we obtain the following theorem.

\proclaim {Theorem 8} For $n \in \Bbb Z_+$, we have
$$\align
&E_{n, \chi, q}^{(h, r)} (x) \\ &=\frac{[2]_q^r}{[2]_{q^f}^r}[f]_q^n
\sum_{a_1, \ldots, a_r =0}^{f-1} (-1)^{\sum_{j=1}^r
a_j}(\prod_{j=1}^r \chi(a_j)) q^{\sum_{j=1}^r (h-j+1)a_j}
\zeta_{q^f, r}^{(h)}(-n, \frac{x+\sum_{j=1}^r a_j}{f}).\endalign$$
\endproclaim
From (18), we note that
$$\align
&\frac{1}{\Gamma(s)}\int_0^\infty F_{q, \chi}^{(h, r)}(-t,x)t^{s-1} dt\\
& \quad=[2]_q^r \sum_{m_1, \ldots, m_r=0}^\infty
\frac{q^{\sum_{j=1}^{r} (h-j+1)m_j}(\prod_{j=1}^r
\chi(m_j))(-1)^{m_1 + \cdots + m_r}}{[m_1 + \cdots + m_r+x]_q^s},
\tag20
\endalign$$
where $s \in \Bbb C$, $x \in \Bbb R$ with $x \neq 0, -1, -2,
\ldots$.

From (20), we define the Dirichlet's type multiple $(h,
q)$-$l$-function associated with the generalized multiple $q$-Euler
polynomials attached to $\chi$.

\proclaim {Definition 9} For $s \in \Bbb C$, $x \in \Bbb R$ with $x
\neq 0, -1, -2, \ldots$, we define the Dirichlet's type multiple
$q$-$l$-function as follows :
$$l_q^{(h)} (s, x |\chi) =[2]_q^r \sum_{m_1, \ldots, m_r=0}^\infty
\frac{q^{\sum_{j=1}^{r} (h-j+1)m_j}(\prod_{i=1}^r \chi(m_i))(-1)^{m_1 +
\cdots + m_r}}{[m_1 + \cdots + m_r+x]_q^s}.$$
\endproclaim
Note that $l_q^{(h)} (s, x |\chi)$ is a meromorphic function in
whole complex plane. It is easy to show that
$$\align
& l_q^{(h)} (s, x |\chi) \\&=\frac{[2]_q^r}{[2]_{q^f}^r}\frac{1}{[f]_q^{s}}
\sum_{a_1, \ldots, a_r =0}^{f-1} (-1)^{\sum_{j=1}^r
a_j}(\prod_{j=1}^r \chi(a_j)) q^{\sum_{j=1}^r (h-j+1)a_j}
\zeta_{q^f, r}^{(h)}(s, \frac{x+\sum_{j=1}^r a_j}{f}).
\endalign$$
By (19) and (20), we obtain the following theorem.

\proclaim {Theorem 10}  For $n \in \Bbb Z_+$, we have
$$l_q^{(h)} (-n, x |\chi)=E_{n, \chi, q}^{(h, r)}(x).$$
\endproclaim

Finally, we give the $q$-extension of Barnes' type multiple Euler
polynomials in (2). For $x, a_1, \ldots, a_r \in \Bbb C$ with
positive real part, let us define the Barnes' type mutiple $q$-Euler
polynomials in $\Bbb C$ as follows :
$$\align
&F_q^{(r)}(t, x|  a_1, \ldots, a_r ; b_1, \ldots, b_r) \\
& \quad = [2]_q^r \sum_{m_1, \ldots, m_r=0}^\infty (-1)^{m_1 +
\cdots + m_r} q^{(b_1 +1)m_1 +
\cdots + (b_r +1)m_r} e^{[a_1 m_1 + \cdots +a_r m_r +x]t} \tag 21\\
& \quad = \sum_{n=0}^\infty E_{n, q}^{(r)}(x| a_1, \ldots, a_r ;
b_1, \ldots, b_r) \frac{t^n}{n!},
\endalign$$
where $b_1, \ldots, b_r \in \Bbb Z$. By (21), we see that
$$\align
&E_{n, q}^{(r)}(x| a_1, \ldots, a_r ; b_1, \ldots, b_r)\\
& \quad=\frac{[2]_q^r}{(1-q)^n} \sum_{l=0}^n
\binom{n}{l}\frac{(-1)^l q^{lx}}{(1+q^{la_1 +b_1
+1})\cdots(1+q^{la_r
+b_r +1})} \\
& \quad =[2]_q^r \sum_{m_1, \ldots, m_r=0}^\infty (-1)^{m_1 + \cdots
+ m_r}q^{(b_1 +1)m_1 + \cdots + (b_r +1)m_r}[a_1 m_1 + \cdots +a_r
m_r +x]_q^n.
\endalign$$
From (21), we note that
$$\align
&\frac{1}{\Gamma(s)}\int_0^\infty F_{q}^{(r)}(-t, x|a_1, \ldots, a_r ; b_1, \ldots, b_r)t^{s-1} dt\\
& \quad=[2]_q^r \sum_{m_1, \ldots, m_r=0}^\infty \frac{(-q)^{m_1 +
\cdots + m_r}q^{b_1 m_1 + \cdots + b_r m_r}}{[a_1 m_1 + \cdots + a_r
m_r+x]_q^s}. \tag22
\endalign$$
By (22), we define the Barnes' type multiple $q$-zeta function as
follows :
$$\align
& \zeta_{q, r} (s, x|a_1, \ldots, a_r ; b_1, \ldots, b_r)\\
& \quad =[2]_q^r \sum_{m_1, \ldots, m_r=0}^\infty \frac{(-q)^{m_1 +
\cdots + m_r}q^{b_1 m_1 + \cdots + b_r m_r}}{[a_1 m_1 + \cdots + a_r
m_r+x]_q^s},
\endalign$$
where $s \in \Bbb C$, $x \in \Bbb R$ with $x \neq 0, -1, -2,
\ldots$. By (21), (22) and (23), we obtain the following theorem.

\proclaim {Theorem 11}  For $n \in \Bbb Z_+$, we have
$$\zeta_{q, r} (s, x|a_1, \ldots, a_r ; b_1, \ldots, b_r)=E_{n, q}^{(r)}(x| a_1, \ldots, a_r ; b_1, \ldots, b_r).$$
\endproclaim

Let $\chi$ be the Dirichlet's character with conductor $f \in \Bbb
N$ with $f \equiv 1 \pmod 2$. Then the generalized Barnes' type
multiple $q$-Euler polynomials attached to $\chi$ are defined by
$$\align
& F_{q, \chi}^{(r)}(t,x|a_1, \ldots, a_r ; b_1, \ldots, b_r)\\&
\quad = [2]_q^r \sum_{m_1, \ldots, m_r=0}^\infty (-q)^{m_1 + \cdots
+m_r}q^{b_1 m_1 +
\cdots + b_r m_r}(\prod_{i=1}^r \chi(m_i)) e^{[a_1 m_1 + \cdots +a_r m_r +x]_q t}\tag 24\\
& \quad = \sum_{n=0}^\infty E_{n, \chi, q}^{(r)}(x| a_1, \ldots, a_r
; b_1, \ldots, b_r) \frac{t^n}{n!},
\endalign$$
From (24), we note that
$$\align
&\frac{1}{\Gamma(s)}\int_0^\infty F_{q, \chi}^{(r)}(-t, x|a_1, \ldots, a_r ; b_1, \ldots, b_r)t^{s-1} dt\\
& \quad=[2]_q^r \sum_{m_1, \ldots, m_r=0}^\infty \frac{(-q)^{m_1 +
\cdots + m_r}q^{b_1 m_1 + \cdots + b_r m_r}(\prod_{i=1}^r
\chi(m_i))}{[a_1 m_1 + \cdots + a_r m_r+x]_q^s}. \tag25
\endalign$$
By (25), we can define Barnes' type multiple $q$-$l$-function in
$\Bbb C$. For $s \in \Bbb C$, $x \in \Bbb R$ with $x \neq 0, -1, -2,
\ldots$, let us define the Barnes' type multiple $q$-$l$-function as
follows :
$$\align
& l_q^{(r)} (s, x |a_1, \ldots, a_r ; b_1, \ldots, b_r) \\&=[2]_q^r
\sum_{m_1, \ldots, m_r=0}^\infty \frac{(-q)^{m_1 + \cdots +
m_r}q^{b_1 m_1 + \cdots + b_r m_r}(\prod_{i=1}^r \chi(m_i))}{[a_1
m_1 + \cdots + a_r m_r+x]_q^s}. \tag26
\endalign$$
Note that $l_q^{(r)} (s, x |a_1, \ldots, a_r ; b_1, \ldots, b_r)$ is
a meromorphic function in whole complex $s$-plane. By (24), (25) and
(26), we easily see that
$$l_q^{(r)} (-n, x |a_1, \ldots, a_r ; b_1, \ldots, b_r)=E_{n, \chi, q}^{(r)}(x| a_1, \ldots, a_r
; b_1, \ldots, b_r)$$ for $n \in \Bbb Z_+$, (see [1-18]).

\vskip 20pt

\quad Taekyun Kim

\quad Division of General Education-Mathematics,
 Kwangwoon University,

 \quad Seoul 139-701, S. Korea
 \quad e-mail:\text{ tkkim$\@$kw.ac.kr}

\vskip 10pt
\quad Young-Hee Kim

\quad Division of General Education-Mathematics,

\quad Kwangwoon University,

\quad Seoul 139-701, S. Korea \quad e-mail:\text{ yhkim$\@$kw.ac.kr}

\enddocument